\documentclass[11pt]{article}

\usepackage{amsbsy,amsfonts,epsf}

\newcommand{\eh}{\hspace{.05in}}  
\newcommand{\R}{\mathbb{R}} 
\newcommand{\C}{\mathbb{C}} 
  
\newcommand{\m}{\frac{1}{2}}

\newcommand{\ds}{\displaystyle}  
\newcommand{\BE}{\begin{equation}}  
\newcommand{\EE}{\end{equation}}  
\newcommand{\af}{\alpha}   
\newcommand{\bt}{\beta}
\newcommand{\ld}{\lambda^+}  
\newcommand{\Ld}{\lambda^-}  
\newcommand{\eps}{\varepsilon}  
  
\newcommand{\A}{\mathcal A}
\newcommand{\CC}{\mathcal C}
 
\newcommand{\Hh}{\mathcal H}
\newcommand{\J}{\mathcal J}
\newcommand{\K}{\mathcal K}
\newcommand{\SK}{\mathcal{SK}}
\newcommand{\T}{\mathcal T}
\newcommand{\X}{\mathcal X}
\newcommand{\sfC}{\sf C}

\newcommand{\br}{\boldsymbol{r}}
\newcommand{\brho}{\boldsymbol{\rho}}
\newcommand{\bw}{\boldsymbol{w}}
\newcommand{\bz}{\boldsymbol{z}}
\newcommand{\deh}{\partial}
\newcommand{\Lim}[1]{\lower5pt\hbox{${{\ds\lim}\atop #1}$}} 
\newcommand{\halfline}[1]{\lower-4pt\hbox{${{\ds\longrightarrow}\atop #1}$}}

\begin{document}
\centerline{\Large{Solving period problems for minimal surfaces}}
\centerline{\Large{with the support function}}
\bigskip

\centerline{By {\it Frank Baginski} at Washington DC and {\it Val\'erio Ramos Batista} at St. Andr\'e}
\bigskip
\centerline{\underline{\hspace{3cm}}}
\bigskip\bigskip

{\bf Abstract.} \ In this paper we show how to bypass the usual difficulties in the analysis of elliptic integrals that arise when solving period problems for minimal surfaces. The method consists of replacing period problems with ordinary Sturm-Liouville problems involving the support function. We give a practical application by proving existence of the sheared Scherk-Karcher family of surfaces numerically described by Wei. Moreover, we show that this family is continuous, and both of its limit-members are the singly periodic genus-one helicoid.
\\
\\
\\
\centerline{\bf 1. Introduction}
\\

In the past decades, the Theory of Minimal Surfaces went through a strong development that started with the works of Douglas [7], Rad\'o [24], Huber [10] and Osserman [22]. While the formers [7,24] contributed to the {\it Theory of Boundary Value Problems}, the latters [10,22] found new important general results on {\it Complete Minimal Surfaces}, that led to the construction of further examples by Chen-Gackstatter [5], Costa [4], Hoffman-Meeks [13], Karcher [15-17], Mart\'\i n and Ramos Batista (see [19] and [25-28]). Karcher was author of the most numerous set of such new surfaces due to a {\it reverse construction method} that he himself developed. In [15-17] one finds the first embedded minimal examples with any given number of helicoidal ends, the first doubly periodic ones out of Scherk's family, the existence proof of Schoen's experimental triply periodic surfaces, genus one saddle towers, and many other original results.     
\\

These constructions are essential for the development of the Global Theory of Minimal Surfaces. Perhaps the most striking was the find of {\it the genus one helicoid} by Hoffman-Karcher-Wei [12] that we call GOH. It is known that a complete embedded minimal surface $S$ of finite total curvature in a flat space must have finite topology (see [20] and [22]). However, except for the plane, the converse is true in $\R^3$ if and only if $S$ has more than one end (see [3]). For a long time, the helicoid had been the sole example with one end until the discovery of the GOH. However, since then there have been no further enrichments in this class of surfaces. 
\\

The GOH was first obtained as a numeric limit of graph surfaces called {\it the twisted singly periodic helicoids with handles}. In 2000, Weber [30] proved these surfaces to exist in a smooth one-parameter family $\Hh_k$, $\m<k\le\infty$. The GOH is the limit member $k=\infty$, while $\Hh_1$ is the well-known {\it singly periodic genus one helicoid} (SGOH), previously found by Hoffman-Karcher-Wei [11]. In [8] and [11], the SGOH is shown to be embedded. Subsequently, Weber finally proved embeddedness of the GOH in [30].    
\\
\begin{figure}[ht]
\centerline{
\epsfxsize 10cm
\hspace{-4cm}
\epsfbox{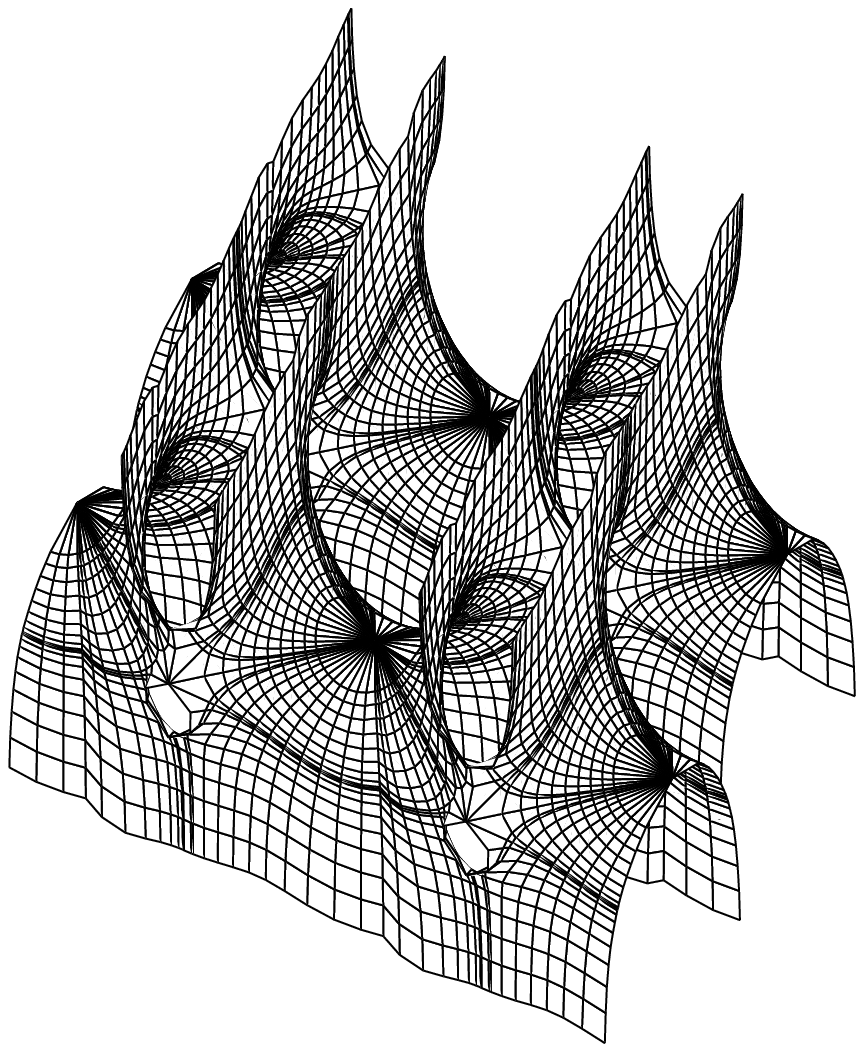}
\epsfxsize 10cm
\hspace{-16.5cm}
\epsfbox{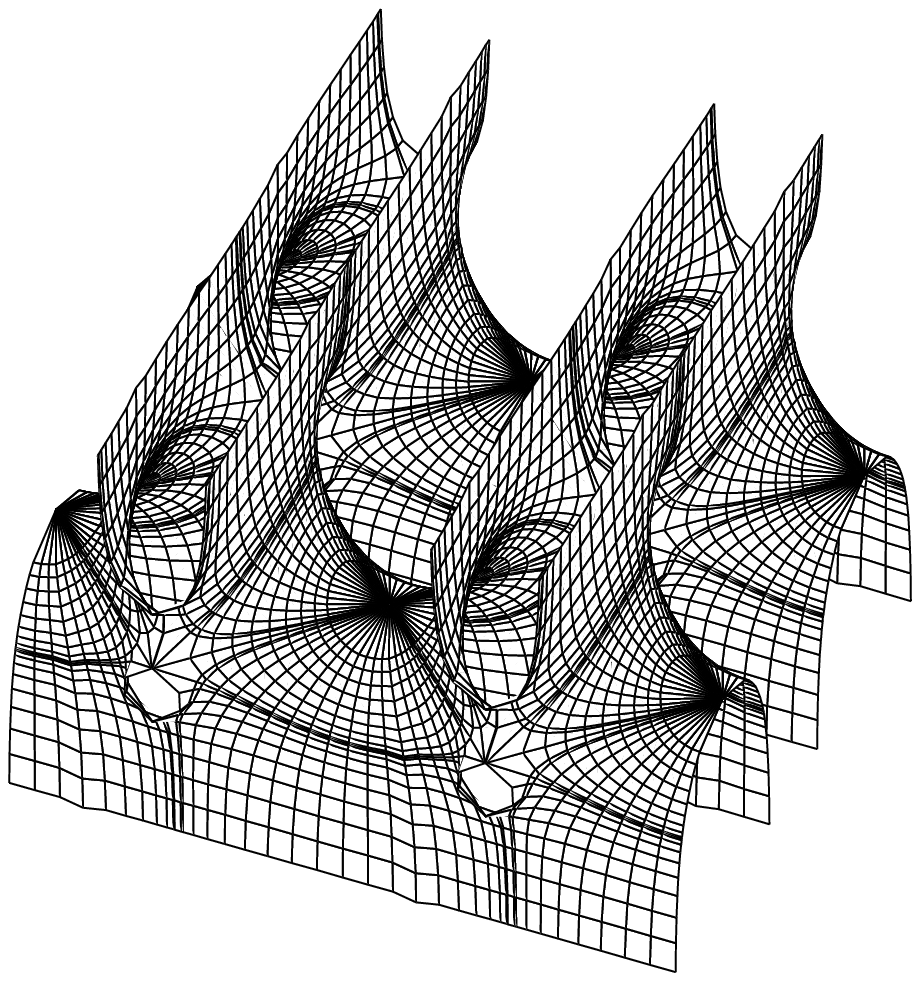}}
\hspace{3cm}(a)\hspace{6cm}(b)
\caption{(a) The Scherk-Karcher surface; (b) one of its sheared deformations.}
\end{figure}

Like the GOH, the SGOH was also first obtained as a numeric limit of graphs. A fundamental fact for its discovery was that Scherk's doubly periodic surface admits sheared deformations which converge to the classical helicoid. A good illustration of this fact can be found in [30]. In the same work the author presents a surface that Karcher obtained by adding to Scherk's surface a handle between every second pair of ends. We call it the Scherk-Karcher surface, which is depicted in Figure 1(a). Using numerical computations, Wei obtained from this example the sheared deformations (see Figure 1(b)), which led to the SGOH. Notice that any sheared Scherk-Karcher surface contains a plane rectangular lattice $\Gamma$ that we can suppose to include the axes $Ox_1$ and $Ox_2$.
\\

Although the limit surface was formally proved to exist, the sheared deformations still remained as purely computational graphs. However, in [14] the authors proved them to exist when {\it close enough} to the standard Scherk-Karcher example. Their periods are difficult to handle and so the embeddedness proof of the SGOH was accomplished by means of other techniques (see [8] or [11]). In fact, this is one of the greatest difficulties at constructing new complete minimal surfaces: the handling of so-called {\it period problems}. In general, one faces transcendental integrals with several interdependent parameters which altogether must fulfil a system of special equalities. If ever solvable, it is usually with extreme difficulties. Therefore, specialists in complete minimal surfaces are turning back to the development of its general theory, looking for non-existence results, or trying new methods of construction.
\\

Like the $\Hh_k$ surfaces, the numeric sheared Scherk-Karcher surfaces are members of a smooth one-parameter family that we call $\SK_k$, starting with Karcher's example at $k=1$ and ending with SGOH at $k=\infty$. In this paper we prove these facts analytically, including embeddedness of the $\SK_k$ surfaces (see Theorem 1.1). We introduce a method which consists of replacing period problems with ordinary Sturm-Liouville problems that are derived from the {\it support function}. This function was first introduced by Minkowski in 1901, conceived as the scalar product $w=\vec X\cdot\vec n$, where $\vec X$ is a parametrisation of a hyper-surface in $\R^m$ with unitary normal $\vec n$. For minimal surfaces in $\R^3$, the support function satisfies the equation 
$$
    \Delta w+2w=0,\leqno(1)
$$
where $\Delta$ is the spherical Laplacian. Equation (1) is a  sufficient local condition for a surface to be minimal. This observation motivated some specialists at the beginning of the 20th century (like Bromwich [2], Darboux [6], and Richmond [29]) to construct new examples by studying some solutions of (1). However, their study was limited to already known minimal surfaces of genus zero, and period problems were not discussed in their works. At that time, period problems were not completely understood. After the works from Hueber [10] and Osserman [22], the minimal surfaces theory took other directions and the support function approach fell into abeyance. Only recently it was used again for an embedded genus-one construction (see [1]).
\\

The herewith presented method is likely to ease period analysis in two different ways. First, one {\it does not} have to solve the Sturm-Liouville Equation. One needs only the {\it sign} of its solution at the final extreme of the definition interval. Second, the coefficients of the ordinary differential equation (ODE) will be real algebraic (or even rational) functions, much easier to handle in comparison with elliptic integrals. Before introducing the method in Section 3, we present a summary of main tools used in this paper in the next section.
\\

The Support Function Method together with some additional geometric analysis will be used to establish the following:
\\

{\bf Theorem 1.1.}\it There exists a continuous one-parameter family of doubly periodic minimal surfaces in $\R^3$ denoted by $\SK_k$, $k\in\R$, having the following properties:

(a) The quotient by its translation group is conformal to a rhombic torus punctured at four points.

(b) For $k\ne 1$, $\SK_k$ contains a rectangular lattice from which one has all the intrinsic symmetries. The member $\SK_1$ coincides with the Scherk-Karcher surface.

(c) $S_k$ is embedded in $\R^3$ and $\SK_{(2-k)}\cong\SK_k$, $\forall$ $k\in\R$. Moreover, $\Lim{k\to\infty}{SK_k}=\,\,$\rm SGOH.
\\

This present work was supported by the following award and grants: NASA NAG5-5353, FAPESP 02/00694-8, FAPESP 07/00569-2 and FAEP 135/05. We are especially grateful to Professor Vera Carrara, University of S\~ao Paulo, for her assistance with some of the differential topology arguments that were utilized in Sections 5 and 6.
\\
\\
\centerline{\bf 2. Background and Notations}
\\

In this section we state some well known theorems on minimal surfaces. For details, we refer the reader to [15], [18], [21] and [22]. In this paper all surfaces are supposed to be regular.
\\

{\bf Theorem 2.1.} (Weierstrass representation). \it Let $R$ be a Riemann surface, $g$ and $dh$ meromorphic function and 1-differential form on $R$, respectively, such that the zeros of $dh$ coincide with the poles and zeros of $g$. Consider the (possibly multi-valued) function $\X:R\to\R^3$ given by
$$
   \X(p):=Re\int^p(\phi_1,\phi_2,\phi_3),\eh\eh where\eh\eh
   (\phi_1,\phi_2,\phi_3):=\m(g^{-1}-g,ig^{-1}+ig,2)dh.\leqno(2)
$$
Then $\X$ is a conformal minimal immersion. Conversely, every conformal minimal immersion $\X:R\to\R^3$ can be expressed like (2) for some meromorphic function $g$ and 1-form $dh$.\rm
\\

{\bf Definition 2.1.} The pair $(g,dh)$ is the \it Weierstrass data \rm and $\phi_1$, $\phi_2$, $\phi_3$ are the \it Weierstrass forms \rm on $R$ of the minimal immersion $\X:R\to\X(R)=S\subset\R^3$.
\\

{\bf Definition 2.2.} A complete, orientable minimal surface $S$ is \it algebraic \rm if it admits a Weierstrass representation such that $R=\bar{R}\setminus\{p_1,\dots,p_r\}$, were $\bar R$ is compact, and both $g$ and $dh$ extend meromorphically to $\bar R$.
\\

{\bf Definition 2.3.} An \it end \rm of $S$ is the image of a punctured neighbourhood $V_p$ of a point $p\in\{p_1,\dots,p_r\}$ such that $(\{p_1,\dots,p_r\}\setminus\{p\})\cap\bar{V}_p=\emptyset$. The end is \it embedded \rm if this image is embedded for a sufficiently small neighbourhood of $p$.
\\

{\bf Theorem 2.2.} \it Let $S$ be a complete minimal surface in $\R^3$. Then $S$ is algebraic if and only if it can be obtained from a piece $\tilde{S}$ of finite total curvature by applying a finitely generated translation group $G$ of $\R^3$.\rm
\\

From now on we consider only algebraic surfaces. The function $g$ is the stereographic projection of the Gau\ss \ map $N:R\to S^2$ of the minimal immersion $\X$. This minimal immersion is well defined in $\R^3/G$, but allowed to be a multivalued function in $\R^3$. The function $g$ is a covering map of $\hat\C$ and the total curvature of $\tilde{S}$ is $-4\pi$deg$(g)$.
\\
\\
\\
\centerline{\bf 3. The Support Function Method}
\\

As explained in the Introduction, if $R$ is a compact Riemann surface punctured at some points and $(g,dh)$ is a Weierstrass pair on it, the corresponding minimal immersion $\X:R\to\R^3$ can take closed curves in $R$ to open curves in $\R^3$. We consider $\X$ as a {\it multivalued} function and $S$ will be invariant under the action of a translation group $G$ in $\R^3$. We are most interested in the case when $\R^3/G$ is still a flat three-dimensional space and $\X:R\to S/G$ is a homeomorphism.
\\

Next, we motivate the {\it support function method}. Consider an analytic regular curve $c:\R\to R$ with $c(t+2)=c(t)$, for any real $t$. In Figure 2 we represented the image under $\X$ of $c|_{[0,2]}$ projected on a plane to be specified later. Suppose this image space curve is invariant under $\varrho$, a $180^\circ$-rotation around a symmetry axis marked with a dashed line. We take $Ox_3$ as this symmetry axis, without loss of generality. Moreover, consider that $\X(1)\in\X([0,2])\cap Ox_3$ and $t=1$ marks exactly the midpoint of the curve. If ``prime'' denotes the derivative with respect to $t$, then $\X'(1)$ is parallel to the plane $Ox_1x_2$. We consider that $\X(0)\in Ox_1x_2$, again without loss of generality. Because of $\varrho$ we have $\X(2)=-\X(0)$, and since $c$ is a periodic function, $\X'(0)=\X'(2)$ with zero third-coordinate. 
\\

\begin{figure}[ht]
\centerline{
\epsfxsize 5cm
\epsfbox{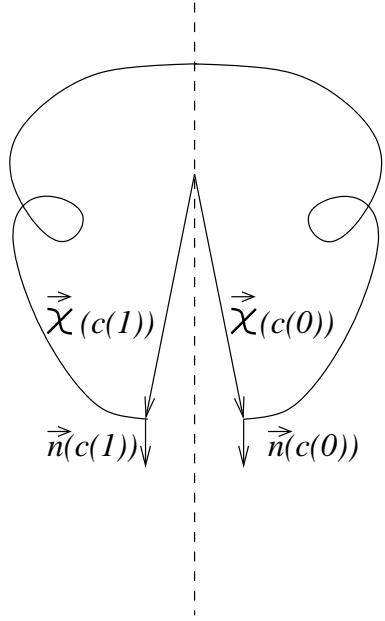}}
\caption{An open period for the minimal immersion $\vec\X$.}
\end{figure}

We also suppose that $Ox_1x_2$ is tangent to $S$ at $\X(0)$, and that $\varrho$ extends to a symmetry of $S$. This happens, for instance, when $S$ is invariant under reflection in $Ox_2x_3$ and $\X\circ c$ is in the plane $Ox_1x_3$. 
\\

From this point on we are going to work with the shifted immersion $X:=\X-\X(0)$, and $w=X\cdot n$. The equalities $\X(t+2)=\X(t)+\X(2)-\X(0)$ and $\X(2)=-\X(0)$, together with the symmetry $\varrho$, imply that $w(t):=w(c(t))$ is even. Hence $w(0)=w'(0)=0$ and $X\circ c$ will close up providing $w$ has a local extreme at $t=1$, and consequently $w'(1)=0$ (later on we shall present a {\it sufficient} condition on $w$ for the periods to close). A simple reckoning gives us
$$
   w''(t)=\vec X(t)\cdot\vec n''(t)-\vec X''(t)\cdot\vec n(t).\leqno(3)
$$

While $X(t)$ is given by a vector of integrals, $X''$ is a purely algebraic expression involving $g(t)$ and $dh(t)$. Let us call $c_{\alpha\bt}$ the coefficients of the first fundamental form in the spherical parametrisation $(u^1,u^2)$ (see [1] for details). Hence
$$
   \vec X = w\vec n + \nabla w\eh\eh{\rm and}\eh\eh n''(t)=-
   \biggl[c_{11}\cdot\biggl(\frac{du^1}{dt}\biggl)^2+c_{22}\cdot\biggl(\frac{du^2}{dt}\biggl)^2\biggl]\vec n+
   \frac{d^2u^1}{dt^2}\vec n_1+\frac{d^2u^2}{dt^2}\vec n_2.\leqno(4)$$

Moreover,
$$
   w\vec n\cdot\vec n''=-w\biggl[
   c_{11}\biggl(\frac{du^1}{dt}\biggl)^2+
   c_{22}\biggl(\frac{du^2}{dt}\biggl)^2
   \biggl]\leqno(5)
$$
and
$$
   \nabla w\cdot\vec n''=
   \frac{\deh w}{\deh u^1}\frac{d^2 u^1}{dt^2}+
   \frac{\deh w}{\deh u^2}\frac{d^2 u^2}{dt^2}.\leqno(6)
$$

If we can get a parametrisation $(u^1(t),u^2(t))$ such that
$$ 
   u^2\equiv{\eh\rm constant\eh\eh and}\eh\eh\frac{du^1}{dt}\ne 0,\leqno(7)
$$
then (6) becomes
$$
   \nabla w\cdot\vec n''=\frac{\deh w}{\deh u^1}\frac{d^2 u^1}{dt^2}=
   \biggl[\frac{\deh w}{\deh u^2}\frac{du^1}{dt}+
   \underbrace{\frac{\deh w}{\deh u^2}\frac{du^2}{dt}}_{=0}\biggl]
   \cdot\biggl(\frac{du^1}{dt}\biggl)^{-1}\cdot\frac{d^2u^1}{dt^2}=
   w'(t)\cdot\biggl[\ln\biggl|\frac{du^1}{dt}\biggl|\biggl]'.\leqno(8)
$$

In this case, (3) together with (5) and (8) will give us the following ODE:
$$
\left\{ \begin{array}{ll}
                w''=p+qw+rw',  \\
                w(0)=w'(0)=0. \\
        \end{array}
        \right.\leqno(9)
$$

Here we have
\[
   p(t)=-\vec X''(t)\cdot\vec n(t)=\vec X'(t)\cdot\vec n'(t);
\]
\[
   q(t)=-c_{11}\biggl(\frac{du^1}{dt}\biggl)^2;\eh{\rm and}
\]
\[
   r(t)=\biggl[\ln\biggl|\frac{du^1}{dt}\biggl|\biggl]'.
\]

As explained in the Introduction, one does not have to find the solution of (9), since we know it {\it is} $w=X\cdot n$. The further condition $w'(1)=0$ holds for a {\it closed} curve. Later on we shall apply the following result from [23, p. 6]:
\\

{\bf Theorem 3.1.} \it Suppose that $w$ is a non-constant solution of the differential inequality $\ddot{w}+G(s)\dot{w}+H(s)w\ge 0$, $s\in[0,1]$, having one-sided derivatives at $0$ and $1$. Suppose further that $H(s)\le 0$ and $\max\{|G|,|H|\}$ is bounded on every closed subinterval of $(0,1)$. In this case, $w$ takes its maximum at $\deh[0,1]$. Moreover, if $w$ has a nonnegative maximum at $1$ and $G(s)-(1-s)H(s)$ is bounded from above at $s=1$, then $\dot{w}(1)>0$.\rm
\\

Now consider $\R^3$ generated by an orthonormal basis $\{e_1,e_2,e_3\}$. In the case of {\it spherical coordinates}, $u^1$ is the angle $\theta$ that the position vector $X$ makes with $e_3$ and $u^2$ is the angle $\varphi$ that $X-(X\cdot e_3)e_3$ makes with $e_1$. For instance, if the curve represented in Figure 2 is an inverse image by $g$ of a meridian in $S^2$, then we can always take $(u^1,u^2)=(\theta,\varphi)$ in such a way that both $\theta(t)$ and $\varphi(t)$ fulfil (7). This makes the right-hand side of (6) equal $r(t)w'(t)$. However, the choice of the $t$-functions $(\theta,\varphi)$ is tied to the parametrisation 
$$
   \vec n=(\sin\theta\cos\varphi,\sin\theta\sin\varphi,\cos\theta)=
   \biggl(\frac{2Re\{g\}}{|g|^2+1},\frac{2Im\{g\}}{|g|^2+1},\frac{|g|^2-1}{|g|^2+1}\biggl).\leqno(10)
$$

From (10) one easily reads off the relation 
$$
   g=\frac{\sin\theta e^{i\varphi}}{1-\cos\theta},\leqno(11)
$$
which {\it must} be used in order to establish $(\theta(t),\varphi(t))$ in terms of $g(t)$. From (11) we have:
$$
   \theta(t)=2\arctan(|g(t)|^{-1})\eh\eh{\rm and}\eh\eh\varphi(t)=\arctan\frac{Im(g)}{Re(g)}.\leqno(12)
$$

On $c(t)$ we write $dh(t)=h'(t)dt$. After rotation about $Ox_3$, $S$ is positioned so that either $\theta\equiv\pi/2$ or $\varphi\equiv 0$. In the former case $|g|\equiv 1$ and $p(t)=\varphi'Re\{ih'\}$, while in the latter $g$ is real and $p(t)=2(\ln|g|)'Re\{h'\}$. As explained at the Introduction, the functions $p$, $q$ and $r$ will be algebraic, or even rational. This is because $g$ and $dh$ can always be written as rational expressions involving two meromorphic functions, which give an algebraic equation for the compact Riemann surface $R$.
\\

At this point we have an important remark. The condition $w'(1)=0$ is necessary, but not sufficient for $X\circ c$ to be closed. Firstly, $t=1$ could be an inflection point of $w$. We shall see that it is indeed an extreme if the Gaussian curvature $K(c(1))\ne 0$, and this will imply that $X(1)$ is parallel to $Ox_2x_3$. Secondly, even in this case $X(1)$ can be non-vertical, and so we have just guaranteed that $X\circ c$ has no periods in the $x_1$-direction. But, for instance, if $X\circ c$ is in the plane $Ox_1x_3$, then $X\circ c$ will already have no periods in the $x_2$-direction. 
\\

The Gaussian curvature for minimal surfaces,
\[
   K=-\frac{16|dg/dh|^2}{|g|^2(|g|+1/|g|)^4},
\]
is strictly negative at a point if and only if $g$ is injective at the point. Suppose $g$ is injective at $c(1)$, hence $\theta'\ne 0$ from (12). We recall that $dh$ is zero whenever $g\in\{0,\infty\}$ at a regular point of the surface. From (11), observe that $t\to 1$ implies $\theta\to 0$. A simple reckoning gives $dw/d\theta|_{\theta=0}=X|_{\theta=0}\cdot (1,0,0)$, namely the first coordinate of the position vector. Therefore, $w'(1)\ne 0$ if and only if $X(1)$ has non-zero first coordinate. For $w'(1)=0$, near $\theta=0$ we can write $w(\theta)=x_3(\theta)+{\cal O}(\theta^2)$. Since $x_3$ has a local extreme at $\theta=0$ because of $\varrho$, then $t=1$ is an extreme for $w$.
\\

Before we conclude this section, it is important to notice that the Support Function Method involves more than the analysis of the boundary value problem defined in (9). It also includes the geometric arguments presented herein, and many times it is possible to deduce the sign of $w'(1)$ even without looking at (9). This is especially true when $X\circ c$ is convex, for suppose that $g$ is real along the curve in Figure 2. Hence $w$ can be computed through its projection onto $Ox_1x_3$. We locally re-parametrise it at the bottom right-hand apex by $(t,f(t))$, for a convex function $f$. The apex is then attained at $(t_0,f(t_0))$, for a certain positive $t_0$. Therefore, $w(t)=(-tf'+f)/\sqrt{1+f'^2}$ and consequently $w'(t_0)=-t_0f''(t_0)<0$. 
\\

Figure 2 shows the case of a ``too long'' period. For a ``too short'' period we would get $w'$ positive at the corresponding apex. These properties will be frequently used in Section 5.
\\
\\
\\
\centerline{\bf 4. Preliminaries}
\\

As explained in the Introduction, the Scherk-Karcher surfaces admit numerical sheared deformations, of which the formal existence was never proved before. This section is devoted to their Weierstrass data, which are obtained by {\it the reverse construction method} from Karcher. There is a large literature about Karcher's method (see, e.g., [15-17], [19] and [25-28]), and so details will be omitted here. The Weierstrass data will be chosen with the help of Figure 1(b), which represents the sought after surface $S$. Figure 3(a) reproduces the quotient of $S$ by its translation group with a shaded fundamental domain. 
\\ 

\begin{figure} [ht] 
\centerline{ 
\epsfxsize 9cm 
\hspace{-1.6in}
\epsfbox{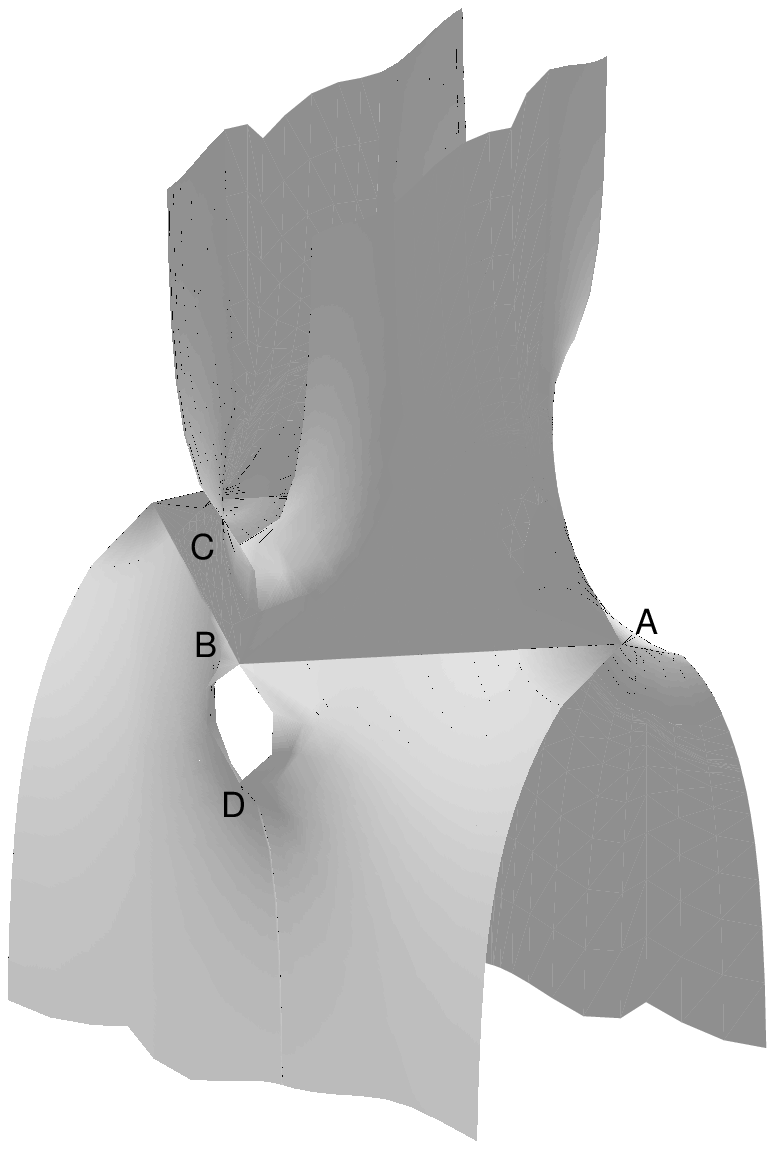}
\epsfxsize 7cm 
\hspace{0.0in}
\epsfbox{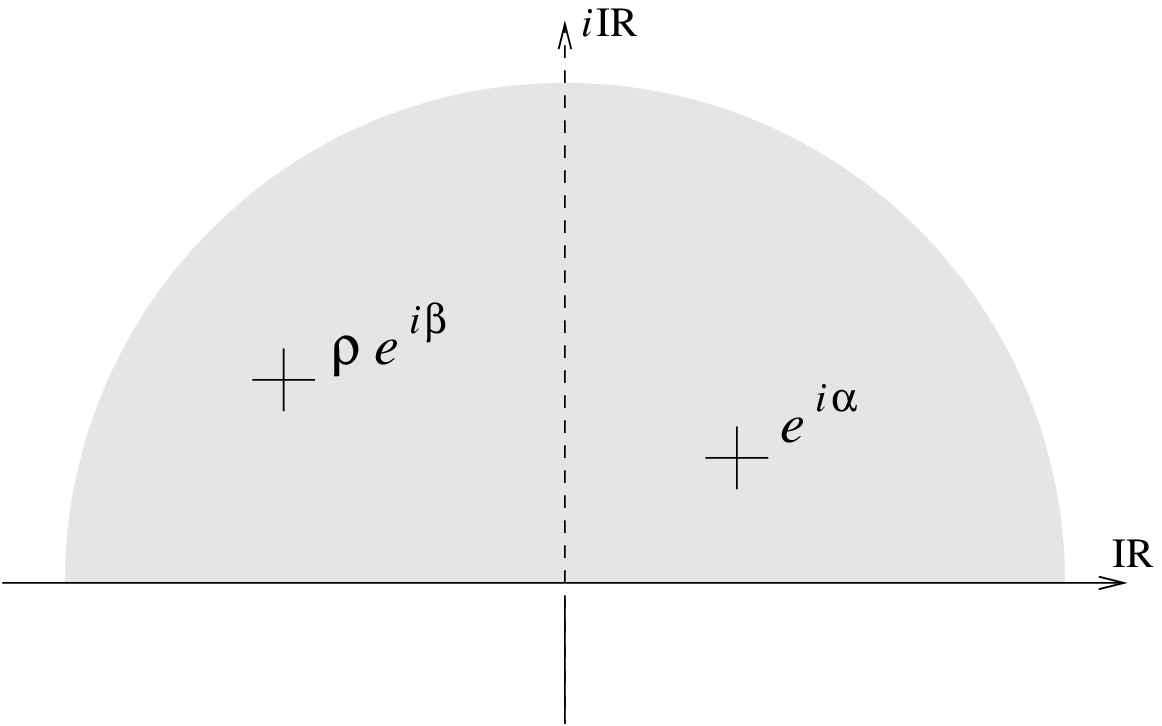}}
\hspace{1in}(a)\hspace{2.5in}(b)
\caption{(a) A fundamental piece of $S$; (b) image values of $z$.}  
\end{figure} 

Notice that Figure 1(b) suggests exactly two periods. We recall the lattice $\Gamma$ defined in the Introduction. The only intrinsic symmetries of $S$ are given by $180^\circ$ rotations around the lines of $\Gamma$, which includes the axes $x_1$ and $x_2$. Therefore, $S$ is invariant under 180$^\circ$ rotations around the vertical axis $x_3$. This rotation has exactly the points $A$, $B$, $C$ and $D$ as fixed points in the quotient of $S$ by its translation group. Together with a compactification of its ends, this quotient is a rhombic torus $T$. Hence we consider the hyperelliptic function $z:T\to\hat\C$ with $z(A)=\infty$, $z(B)=0$, $z(C)=e^{i\af}$ and $z(D)=e^{-i\af}$, where $\af\in(0,\pi)$. An algebraic equation of $T$ is given by $u^2=z+1/z-2\cos\af$.
\\

Figure 4 shows a fundamental domain of the torus lattice in $\C$. We remark that this lattice is different from $\Gamma$, since $T$ is rhombic. However, the bold straight lines in Figure 4 indicate the corresponding fundamental rectangle of $\Gamma$. The letters with prime are the ones omitted in Figure 3(a).\\

\begin{figure} [ht] 
\centerline{ 
\epsfxsize 10cm 
\epsfbox{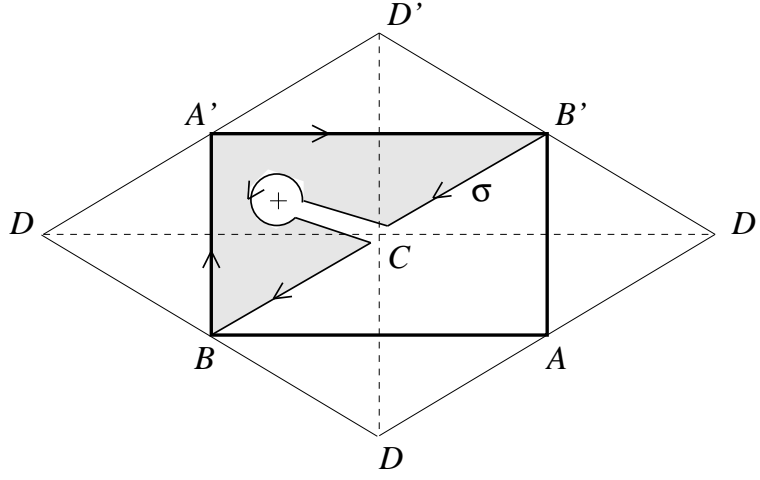}}
\caption{The torus $T$ with the curve $\sigma$.}  
\end{figure} 

We now list some important properties of $z$, which are verified in [11] or [12]. Exactly on $B\to A'$ and $B'\to A$, $z$ is real and negative, while it is positive on $A\to B$, $A'\to B'$ and nowhere else. Exactly on the dashed lines, $z$ is unitary. Therefore, $z=-1$ at the intersection between $B\to A'$ and $\overline{CD}$, while $z=1$ at $\overline{AB}\cap\overline{CD}$.
\\

We recall that the Weierstrass data $(g,dh)$ are such that $g$ is the stereographic projection of the unitary normal on $S$. In Figure 3(a), $B$ is the origin, $Ox_2=\halfline{\ds BA}$ and $Ox_1$ is left-to-right oriented. Therefore, based on Figure 3 we obtain the following relation:  
$$    
   g=u/\sqrt{\kappa}, {\rm \ for \ some \ positive \ constant \ } \kappa.\leqno(13)
$$ 

From (13) it is very easy to verify that $g$ is consistent with the behaviour of the unitary normal on the expected symmetry lines. Now we define $dh$ according to the kind of ends and regular points of $S$. The ends of $S$ are attained at $z=\rho e^{\pm i\bt}$. From Figure 3, by inspection we obtain
$$
   dh=\frac{-idz/z}{z/\rho+\rho/z-2\cos\bt}.\leqno(14)
$$

Now we have (13) and (14), namely Weierstrass data $(g,dh)$ on $T$. Let us define $\T:=T\setminus z^{-1}(\{\rho e^{\pm i\bt}\})$ and consider the corresponding minimal immersion $X:\T\to\R^3$. Well known arguments (see [15-17] or [25-28]) make us conclude that $S=X(\T)$ has all of the expected symmetry lines, Scherk-ends and regularity, according to what Figures 1(b) and 3(a) suggest. However, it is not true that $X(\T)$ has only two periods, unless the parameters $\af$, $\bt$, $\kappa$ and $\rho$ assume some suitable values. In fact, if $\gamma$ is a curve connecting $0$ and $e^{i\af}$ in the upper half-plane, then one must have
$$
   Re\int_\gamma(\phi_1,\phi_2)=0.\leqno(15)
$$

A simple reckoning shows that (15) is equivalent to
$$
   \kappa=\frac{\int_\gamma udh}{\int_\gamma\overline{dh/u}}.\leqno(16)
$$      

However, the ends of $S$ are not supposed to intersect. Namely, $g$ is unitary for $z=\rho e^{\pm i\bt}$. If $z=\rho e^{i\bt}$, then $u=\pm\kappa^{1/2}e^{i\omega}$, $\omega\in\R$, where 
$$
   \kappa e^{2i\omega}=(\rho+1/\rho)\cos\bt-2\cos\af+i(\rho-1/\rho)\sin\bt.\leqno(17)
$$

This again determines $\kappa$ as another function of $\af$, $\bt$ and $\rho$. Therefore (16) becomes
\begin{eqnarray*}
\kappa & = & \{[(\rho+1/\rho)\cos\bt-2\cos\af]^2+[(\rho-1/\rho)\sin\bt]^2\}^\m \\
       & = & \ds{\int_\gamma\frac{i(z+1/z-2\cos\af)^\m dz}{(z/\rho+\rho/z-2\cos\bt)z}\biggl/\int_\gamma\overline{\frac{i(z+1/z-2\cos\af)^{-\m}dz}{(z/\rho+\rho/z-2\cos\bt)z}}}.\\
\end{eqnarray*}

Of course, for $z=\rho e^{-i\bt}$ we have $g=\pm e^{-i\omega}$, as expected due to the straight lines of the surface. However, $g=\pm e^{\pm i\omega}$ at the ends still does not guarantee that $\big|\tan|\omega|\big|$ equals the ratio between the corresponding sides of the rectangle, as Figure 3(a) suggests. For this to be true, (15) is sufficient, as we shall prove now.
\\

In Figure 4, the point marked with a ``+'' represents one of the four Scherk-ends given by $(z,u)\in\{(\rho e^{i\bt},\pm\kappa^{1/2}e^{i\omega}),(\rho e^{-i\bt},\pm\kappa^{1/2}e^{-i\omega})\}$. An easy computation shows that
$$
   2\pi i Res(\phi_1,(z,u)=(\rho e^{i\bt},\kappa^{1/2}e^{i\omega}))=-\pi\sin\omega\cdot\csc\bt;\leqno(18)
$$
$$
   2\pi i Res(\phi_2,(z,u)=(\rho e^{i\bt},\kappa^{1/2}e^{i\omega}))=\pi\cos\omega\cdot\csc\bt;\leqno(19)
$$ 
$$
   2\pi i Res(\phi_3,(z,u)=(\rho e^{i\bt},\kappa^{1/2}e^{i\omega}))=i\pi\csc\bt.\leqno(20)
$$

Consider $\sigma$ to be the closed oriented curve containing stretches $B\to A'$ and $A'\to B'$, as illustrated in Figure 4. Notice that $\phi_1$ is first real and then purely imaginary on $B\to A'$ and $A'\to B'$, respectively. However, $\phi_2$ is first purely imaginary and then real on $B\to A'$ and $A'\to B'$, respectively. Since $\sigma$ is null-homotopic, the above residues will match the segment lengths if (15) holds.  
\\

In the next section one solves period problems by means of the support function.
\\
\\
\\
\centerline{\bf 5. Practical Application of the Support Function Method}
\\

As we have mentioned already, the periods will close up if (15) holds for a curve $\gamma$ connecting $0$ and $e^{i\af}$ in the upper half-plane. This curve can be the stretch $B\to C$ of $\sigma$ depicted in Figure 4. Now consider $D'\to B'$. The involution on $T$ which fixes $\overline{AB}$ is given by $(z,u)\to(\bar{z},\bar{u})$. This involution corresponds to $180^\circ$ rotation around the segments of $S$ which are parallel to $Ox_2$. Therefore, the integral of $\phi_1$ on either $D'\to B'$ or $B'\to C$ is the same. The path $D'\to B'\to C$ will be homotopically equivalent to $D'\to C$ if no Scherk-end is contained in the triangle $D'B'C$. Otherwise, if in the triangle we have $(\rho e^{i\bt},\kappa^{1/2}e^{i\omega})$, there is also $(\rho e^{-i\bt},\kappa^{1/2}e^{-i\omega})$ inside, and no Scherk-end is contained in the {\it other} triangle $DBC$ with horizontal $D\to C$.
\\  

None of these triangles will contain Scherk-ends if $\rho>1$. In fact, the standard Scherk-Karcher surface exists for $\af=\bt=\pi/2$ and a {\it unique} $\rho_1>1$, as proved in Appendix A. We remark that neither [12] nor [14] establishes this result. Since we are looking for a continuous family $\SK_k$, $1\le k<\infty$, starting at $\rho=\rho_1$, the restriction $\rho\ge 1$ in our Weierstrass data (13) and (14) can be imposed without harm. Indeed, at the end of the previous section we saw that (15) guarantees the $\Gamma$-ratio to be $\big|\tan|\omega|\big|$. None of the $\SK$-surfaces could have $\rho=1$, otherwise by (17) it followed $\omega=0$ mod $\pi/2$ and $\Gamma$ would degenerate. Later we shall see that $\rho_k$ converges to a finite value in $(1,\infty)$ when $k$ diverges to $\infty$. From now on consider $\rho\in[1,\infty)$. 
\\

We shall then first take $\rho>1$ and the segment $D\to C$ upwards, on which $z(t)=e^{it}$, $-\af\le t\le\af$. Consider the curve $c:\R\to T$ such that $c|_{[-1,3]}$ is one of the generators for $H_1(T)$, namely $D\to C\to D'$. A simple reckoning shows that $X$ induces the {\it hyperelliptic involution} $(z,u)\to(z,-u)$ on $S$. This means that $S$ is invariant under $180^\circ$ rotation around a vertical axis through any of the image points $A$, $B$, $C$ or $D$. Now, according to Section 3 we should position the extremes of $X\circ c$ on $Ox_1x_2$. However, $S$ is invariant by $180^\circ$ rotation around any of its segments parallel to $Ox_2$. Therefore, up to re-parametrisation we have $c(0)=\overline{DC}\cap\overline{AB}$ and $w(0)=w'(0)=0$. Notice that $z(c(0))=1$ and $X(0)$ is the origin.    
\\

At this point, we are ready to consider the ODE (9). From (13) and (14) one easily reads off the functions $p$, $q$ and $r$ for $z(t)=e^{it}$, $0\le t\le\af$:
\[
   p(t)=\frac{-\sin t}{\cos t-\cos\af}
   \cdot\frac{(1/\rho+\rho)\cos t-2\cos\bt}{[(1/\rho+\rho)\cos t-2\cos\bt]^2+[(1/\rho-\rho)\sin t]^2};
\]
\[
   q(t)=-\frac{2\kappa\sin^2 t/(\cos t-\cos\af)}
   {[\kappa+2(\cos t-\cos\af)]^2};
\]
\[
   r(t)=\cot t+\frac{0.5\sin t}{\cos t-\cos\af}+\frac{2\sin t}{\kappa+2(\cos t-\cos\af)}.
\]

Notice that $r$ has a singularity at $t=0$, but since $w'$ is real analytic with $w'(0)=0$, then $rw'$ is finite at $t=0$. Now take the extreme value $\rho=1$. If $\bt<\af$, then $p$ will have a singularity at $t=\bt$. However, this case can be treated geometrically. Since $\rho=1$, there are additional symmetries which come out on the surface. The image curves of $z(t)=e^{it}$, $0\le t<\bt$ and $\bt<t\le\af$, are both parallel to the plane $x_2=0$. For $\bt<t\le\af$ we have a convex curve, symmetric with respect to a plane parallel to $x_1=0$ and with vertex at $t=\af$. Therefore $w'(\af)>0$, according to the argument at the end of Section 3. 
\\

At this point we remark that $\af$ could have been restricted to the interval $(0,\pi/2]$. If one took $\pi/2<\af<\pi$, then $\tilde{\af}:=\pi-\af$ together with the involution $(z,u)\to(-z,iu)$ would give us $ig=iu/\sqrt{\kappa}$, $(iu(-z))^2=z+1/z-2\cos\tilde{\af}$. In $\R^3$, it just means a $90^\circ$ rotation of the surface around $Ox_3$. In other words, this would simply change our conventions of $X(\overline{A'B})$ on $Ox_1$ to $Ox_2$ and so on. Therefore, without loss of generality, henceforth we take $\af\in(0,\pi/2]$. 
\\

If $\bt>\af$, it is another plane curve of symmetry which comes out, this time entirely given by $z(t)=e^{it}$, $0\le t\le\af$. The image curve is convex and contained in a plane parallel to $x_2=0$. However, this case cannot be treated just geometrically and so we shall make use of classical analytic arguments. On $z(t)$ one has $Re\int\phi_1=0$ if and only if $I_1=I_2$, where
\[
   I_1:=\int_0^\af\biggl(\frac{\cos t-\cos\af}{\cos\af-\cos\bt}\biggl)^\m\cdot\frac{dt}{\cos t-\cos\bt}
\]
and
\[
   I_2:=\int_0^\af\biggl(\frac{\cos\af-\cos\bt}{\cos t-\cos\af}\biggl)^\m\cdot\frac{dt}{\cos t-\cos\bt}.
\]

On the one hand, for $\bt$ approaching $\af$, $I_1$ diverges to $+\infty$ while $I_2$ remains finite (see Appendix B). On the other hand, by taking into account that $\cos t-\cos\af<1+\cos\af$, for $\bt$ approaching $\pi$ we get $I_1<I_2$. Since $Re\int\phi_1=I_2-I_1$, $\bt\cong\af$ implies $w'(\af)>0$ while $\bt\cong\pi$ implies $w'(\af)<0$. The {\it intermediate value theorem} assures the existence of a certain $\bt^+\in(\af,\pi)$ at which $w'(\af)=0$. From Appendix C, it follows the existence of a curve $\CC^+$ in the region $\A:=\{z\in\C:|z|>1>1-Im\{z\}\}$ with $w'|_{\CC^+}\equiv 0$. Since $w$ is real analytic, this curve separates $\A$ into finitely many simply connected regions. Also in Appendix C, we show that $\CC^+$ cannot diverge to $\infty$.
\\

In order to interpret this fact geometrically, observe that the immersion $X$ takes $\overline{BA'}$ to a segment in $Ox_1$. Moreover, (17) for $\rho=1$ implies that $\omega=\pi/2$ mod $\pi$, and from (18)-(20), the upper Scherk-ends measure $\pi\csc\bt$. If they are shorter than $X(\overline{A'B})$, then $w'(\af)$ is positive. If they are larger, then $w'(\af)<0$. 
\\

From the Support Function Method, if we choose $\rho e^{i\bt}\in\CC^+$ the point $C$ will project perpendicularly on $Ox_2$. Namely, $X(c(\af))$ has zero first-coordinate, as explained in Section 3. Therefore, there is a vertical axis $Ox_3'\perp Ox_2$ such that $C\in Ox_3'$. Hence $180^\circ$ rotation around $Ox_3'$, denoted $\varrho'$, is a symmetry of $S$. Thus, the image of $X(c(0))$ under $\varrho'$ lies again in $Ox_2$. But since $S$ is invariant under rotation around $Ox_2$, the curve $X\circ c$ must be closed. This solves the first period problem.\\

Before going ahead, notice that $\bt^+$ is unique because $(\cos\af-\cos\bt)^\m I_1$ decreases while $(\cos\af-\cos\bt)^\m I_2$ increases with $\bt$.
\\ 

We take now the horizontal segment $D\to C$, on which $z(t)=e^{i(\pi-t)}$, $\af-\pi\le t\le\pi-\af$. Consider the curve $c:\R\to T$ such that $c|_{[-1,3]}$ is the other generator of $H_1(T)$. Since $S$ is invariant by $180^\circ$ rotation around any of its segments parallel to $Ox_1$, up to re-parametrisation we have $c(0)=\overline{DC}\cap\overline{A'B}$ and $w(0)=w'(0)=0$. Notice that this time $z(c(0))=-1$, but $X(0)$ is again the origin.    
\\

One considers now the curve $z(t)=e^{i(\pi-t)}$, $t\in[0,\pi-\af]$. Hence we have 
\[
   p(t)=\frac{-\sin t}{\cos t+\cos\af}
   \cdot\frac{(1/\rho+\rho)\cos t+2\cos\bt}{[(1/\rho+\rho)\cos t+2\cos\bt]^2+[(1/\rho-\rho)\sin t]^2};
\]
\[
   q(t)=-\frac{2\kappa\sin^2 t/(\cos t+\cos\af)}
   {[\kappa+2(\cos t+\cos\af)]^2};
\]
\[
   r(t)=\cot t+\frac{0.5\sin t}{\cos t+\cos\af}+\frac{2\sin t}{\kappa+2(\cos t+\cos\af)}.
\]

Now take again the extreme value $\rho=1$. If $\bt>\af$, then $p$ will have a singularity at $t=\pi-\bt$. Once more we can use purely geometrical arguments. Since $\rho=1$, the additional symmetries are this time given by $z(t)=-e^{-it}$, $0\le t<\pi-\bt$ and $\pi-\bt<t\le\pi-\af$, both plane curves parallel to $x_1=0$. The stretch $\pi-\bt<t\le\pi-\af$ is convex, symmetric with respect to a plane parallel to $x_2=0$ and with vertex at $t=\pi-\af$. Therefore $w'(\pi-\af)>0$. 
\\

If $\bt<\af$, $z(t)$ leads to another curve of symmetry, this time in a plane parallel to $x_1=0$. The immersion $X$ takes $\overline{BA}$ to a segment in $Ox_2$. Now (17) implies that $\omega=0$ mod $\pi$, but from (18)-(20) the upper Scherk-ends measure again $\pi\csc\bt$. The length of $X(\overline{BA})$ is easily computable as 
\[
   \m\int_0^\infty\biggl[\frac{\kappa^\m}{(t+1/t-2\cos\af)^\m}+\frac{(t+1/t-2\cos\af)^\m}{\kappa^\m}\biggl]
   \frac{dt/t}{t+1/t-2\cos\bt}, 
\]
with $\kappa=2(\cos\bt-\cos\af)$ by (17). At this point, recall the arguments used to analyse the previous case $z(t)=e^{it}$ and $\bt>\af$. Back to our present case, where $z(t)=e^{i(\pi-t)}$ and $\bt<\af$, analogous arguments can be applied. We then conclude that the Scherk-ends are {\it shorter} than $X(\overline{BA})$ only when $\bt>\bt^-$, for a certain $\bt^-\in(0,\af)$, while $w'(\pi-\af)<0$ for $\bt<\bt^-$. Therefore, $\bt^-<\bt<\af$ and $\rho=1$ imply $w'(\pi-\af)>0$. Similarly to the previous case in Appendix C, we get $w'(\pi-\af)=-\infty$ for enough large $\rho$ but $\af$ close to $\pi/2$. However, $w'(\pi-\af)=+\infty$ for $\rho>>1$ and $\af\cong 0$. In any case, the {\it intermediate value theorem} assures the existence of a curve $\CC^-$ in $\A$ with $w'|_{\CC^-}\equiv 0$. This curve also separates $\A$ into finitely many simply connected regions.
\\

From the Support Function Method, if we choose $\rho e^{i\bt}\in\CC^-$, then the point $C$ will project perpendicularly on $Ox_1$. Namely, $X(c(\pi-\af))$ has zero second-coordinate, as explained in Section 3. Therefore, there is a vertical axis $Ox_3''\perp Ox_1$ such that $C\in Ox_3''$. Hence $180^\circ$ rotation around $Ox_3''$, denoted $\varrho''$, is a symmetry of $S$. Thus, the image of $X(c(0))$ under $\varrho''$ lies again on $Ox_1$. But since $S$ is invariant under rotation around $Ox_1$, the curve $X\circ c$ must be closed. This solves the second period problem.
\\

We must be careful at this point because one did not verify yet whether $\CC^+$ and $\CC^-$ eventually intersect. In other words, it still lacks a {\it simultaneous} solution for the first and second periods. A priori, it could happen that $w'(\af)$ close to $(\rho,\bt)=(1,\af)$ changes sign, but we are going to show that this is not the case. Let us look at $w'(\af)$ and $w'(\pi-\af)$ as functions of $(\rho,\bt)$. For $\af=\pi/2$, the curves $\CC^+$ and $\CC^-$ will be symmetric by a reflection in $i\R$. From Appendix D, if we prove that $w'(\pi/2)$ is positive in a punctured neighbourhood of $(1,\pi/2)$, this will mean that $\CC^+\cap\CC^-\ni i\rho_1$, for a certain $\rho_1>1$. In fact, this $\rho_1$ will be the unique $\rho$-value that defines the standard Scherk-Karcher surface $\SK_1$.    
\\

Now we show that neither $\CC^+$ nor $\CC^-$ gets close to $(\rho,\bt)=(1,\af)$. The segments $X(\overline{A'B})$ and $X(\overline{BA})$ measure 
$$
   \m\int_0^\infty\biggl[\frac{\sqrt{\kappa}}{\sqrt{t+1/t\pm 2\cos\af}}+\frac{\sqrt{t+1/t\pm 2\cos\af}}{\sqrt{\kappa}}\biggl]
   \frac{dt/t}{t/\rho+\rho/t\pm 2\cos\bt},\leqno(21)
$$
respectively. Therefore, in a punctured neighbourhood $V$ of $(\rho,\bt)=(1,\af)$, the integrals at (21) diverge to $+\infty$, while the length of the Scherk-ends remain bounded by (18) and (19). Along $z(t)=e^{it}$, $0\le t\le\af$, the function $g$ is real and varies monotonically. Therefore, the projection of $X\circ z$ on the plane $x_2=0$ is convex. Moreover, $w=X\cdot n=({\rm proj}_{x_2=0}X)\cdot n$ and so, according to the argument at the end of Section 3, $w'(\af)>0$ on $V$. Similarly, on $z(t)=-e^{-it}$, $0\le t\le\pi-\af$, $g$ is purely imaginary and varies monotonically. This time ${\rm proj}_{x_1=0}(X\circ z)$ is convex and the corresponding $w$ on this curve has a positive derivative at $t=\af$, for any $(\rho,\bt)\in V\setminus\{(1,\af)\}$.
\\

From the previous arguments we infer that $\CC^+\cap\CC^-\ni i\rho_1$ in the case $\af=\pi/2$, for a certain $\rho_1>1$. In this case $\bt=\pi/2$ and equations (21) turn out to be the same, so the periods close up if and only if (21) equals the absolute value of (18) or (19), namely $\pi/\sqrt{2}$. In other words, the following equality must hold for $\kappa=\rho-1/\rho$:
$$
   \pi\sqrt{2}=\int_0^\infty
   \biggl[\frac{\sqrt{\kappa}}{\sqrt{t+1/t}}+\frac{\sqrt{t+1/t}}{\sqrt{\kappa}}\biggl]
   \frac{dt/t}{t/\rho+\rho/t}.\leqno(22)
$$
   
From Appendix A, it follows that (22) holds for a single $\rho_1>1$. Now denote $S^1_+:=\{z\in S^1: Im\{z\}\ge 0\}$. On the set $(-\infty,-1]\cup S^1_+\cup[1,+\infty)$, the extremes of $\CC^+$ and $\CC^-$ lie alternately. Since they are real analytic curves, their intersection consists of a finite number points. Moreover, from [9, p. 132] their {\it intersection number} is always equal to one (the total summation after attributing sign and degree to each crossing and tangent point). Indeed, by considering $\CC^+$ as an immersed submanifold of $S^2$, we can smoothly join its extremes with a simple curve in $S^2\setminus\A$. Afterwards, one eliminates self-intersections according to the procedure described in [9, p. 127]. The same can be done to $\CC^-$. One gets a family of $S^1$-embeddings in $S^2$, transversal up to arbitrarily small perturbations. From [9, p. 132], their intersection number is zero, a topological invariant. By deducting the single crossing at $S^2\setminus\A$, one gets $\sharp(\CC^+,\CC^-)=1$.       
\\

We start at $\af=\pi/2$ and let $\af$ converge to zero. Figure 5 shows a possible failure at trying to get a continuous family of surfaces parametrised by $\af$.
\\

\begin{figure} [ht] 
\centerline{ 
\epsfxsize 12cm 
\epsfbox{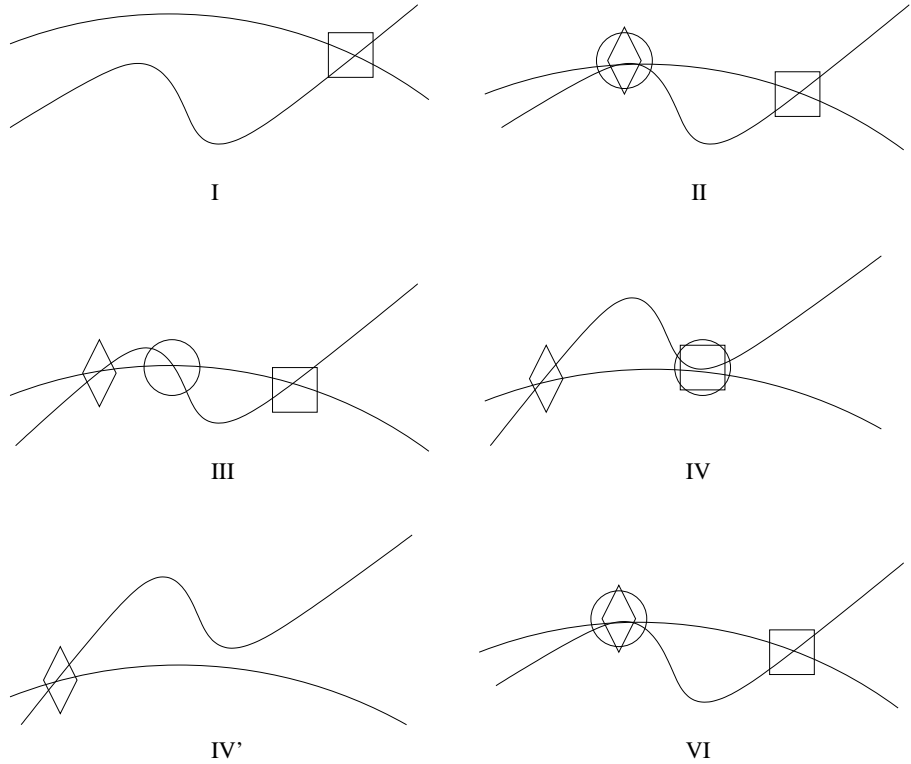}}
\caption{Possible sequence of intersections between $\CC^+$ and $\CC^-$.}  
\end{figure} 

From Figure 5 one sees that the $\Box$-crossing will die off after step IV. However, if we take back steps V:=III and VI:=II by tracking the $\bigcirc$-crossing, we can re-start at VI with the $\Diamond$-crossing. From that step on we take VII:=III and VIII:=IV. In the next section we shall formalise this procedure. It illustrates the fact that one still gets a continuous family of surfaces, now parametrised by a variable that we call $k\in[1,\infty)$. One starts at $\af(1)=\pi/2$ and ends with $\af(\infty)\in(0,\pi/2)$. If $\CC^+\cap\CC^-$ is always transversal, then $\af$ will be a monotone function of $k$, which seems to hold numerically. However, this fact is far from being trivial to prove. Anyway, in the next section we demonstrate that $k\to\infty$ will lead to the SGOH surface, the {\it singly periodic genus one helicoid}.    
\\

Here we summarise what was obtained so far: For every $\af\in(0,\pi/2]$ there are two analytic curves $\CC^{\pm}\subset\A$, along which one closes up each of the two periods. The curves have alternating endpoints for $\af$ in a neighbourhood of $\pi/2$ and vary analytically with $\af$. From their common intersections we can describe a continuous family of surfaces. At this point we have already proved items (a) and (b) of Theorem 1.1. 
\\
\\
\\
\centerline{\bf 6. Limits and Embeddeness}
\\

In the previous section we obtained a continuous one-parameter family of doubly periodic minimal surfaces $\SK_k$ with $k\in[1,\infty)$. Let us now analyse what happens when $k$ diverges to $\infty$. In the general case, one recalls (17)-(19) and (22) becomes a system of two equations:
$$
   2\pi\frac{\cos\omega}{\sin\bt}=\int_0^\infty
   \biggl[\frac{\sqrt{\kappa}}{\sqrt{t+1/t-2\cos\af}}+
   \frac{\sqrt{t+1/t-2\cos\af}}{\sqrt{\kappa}}\biggl]
   \frac{dt/t}{t/\rho+\rho/t-2\cos\bt}\leqno(23)
$$
and
$$
   2\pi\frac{\sin\omega}{\sin\bt}=\int_0^\infty
   \biggl[\frac{\sqrt{\kappa}}{\sqrt{t+1/t+2\cos\af}}+
   \frac{\sqrt{t+1/t+2\cos\af}}{\sqrt{\kappa}}\biggl]
   \frac{dt/t}{t/\rho+\rho/t+2\cos\bt},\leqno(24)
$$
where
$$
   \omega=\left\{
   \begin{array}{ll}
   \m\arctan\frac{(\rho-1/\rho)\sin\bt}{(\rho+1/\rho)\cos\bt-2\cos\af}, 
   \ {\rm if} \ (\rho^2+1)\cos\bt\sec\af\ge 2\rho;\\
   \frac{\pi}{2}+\m\arctan\frac{(\rho-1/\rho)\sin\bt}{(\rho+1/\rho)\cos\bt-2\cos\af}, 
   \ {\rm if} \ (\rho^2+1)\cos\bt\sec\af<2\rho.
   \end{array}\right.\leqno(25)
$$

Both (23) and (24) will simultaneously hold for a space curve $(\af(k),\bt(k),$ $\rho(k))$, $k\in[1,\infty)$, with $\af(1)=\bt(1)=\pi/2$ and $\rho(1)$ equals a certain $\rho_1>1$. Whilst (23) solves the period problem in the $x_1$-direction, (24) solves it in the $x_2$-direction. We recall that, for any fixed $k$, $\CC^+$ starts at $\bt^+\in(\af,\pi)$. From (23) one sees that $\bt^+$ can never belong to the interval $(0,\pi/2]$, because then $\rho=1$ and $w=\pi/2$, and so $\bt^+\in(\pi/2,\pi)$. At the other extreme of $\CC^+$, from (17) a simple reckoning gives
\[
   \lim_{\bt\to 0}\sqrt{\kappa}\sin\omega\csc\bt=
   \frac{(\rho-1/\rho)/2}{\sqrt{\rho+1/\rho-2\cos\af}}.
\] 

From (24), one sees that $\rho$ cannot be close to 1 when $\CC^+$ approaches the real positive axis. In Section 5 we took the paths $z(t)=e^{it}$, $t\in[0,\af]$, and $z(t)=e^{i(\pi-t)}$, $t\in[0,\pi-\af]$. Let us write the integration of $\phi_1$ on the 1st path as $\int_0^\af\phi_1$, and of $\phi_2$ on the 2nd as $\int_0^{\pi-\af}\phi_2$. The simultaneous solution (23) and (24) are then equivalent to 
$$
   Re\int_0^{\af}\phi_1=0\leqno(26)
$$ 
and
$$
   Re\int_0^{\pi-\af}\phi_2=0\leqno(27)
$$  
respectively. Now suppose that $\CC^-$ will always have points with $\rho\ge 1+\eps$, for a certain positive $\eps$, no matter how close $\af$ is to zero. By explicitly writing down $2Re\int_0^{\af}\phi_1$, we get
$$ 
   \int_0^\af
   \frac{\biggl(\frac{\sqrt{\kappa/2}}{(\cos t-\cos\af)^\m}-
   \frac{(\cos t-\cos\af)^\m}{\sqrt{\kappa/2}}\biggl)\cdot
    [(\rho+1/\rho)\cos t-2\cos\bt]dt}
   {[(\rho+1/\rho)\cos t-2\cos\bt]^2+[(\rho-1/\rho)\sin t]^2}.\leqno(28)
$$

By taking $t=\af s$ in (28), one easily computes
\[
   \liminf_{\af\to 0}2\sqrt{\kappa}Re\int_0^{\af}\phi_1\ge\pi,
\]
which contradicts (26). Therefore, $\CC^-$ degenerates to point $1\in S^1$ when $\af$ approaches 0. By joining our conclusions about $\CC^{\pm}$, we see that there is an $\af^*\in(0,\pi/2)$ for which both curves cross at a point $\rho^*\in(1,\infty)$. With $\bt^*=0$, this will give us the same Weierstrass data of the SGOH, as shown further on. According to [8] the SGOH is {\it unique}, and so will be $(\af^*,\rho^*)$. This means that $\CC^+\cap\CC^-$ has an odd intersection number for $\af\in[\af^*,\pi/2]$, and an even one for $\af\in(0,\af^*)$. We then have $(\af,\bt,\rho)|_{k=\infty}:=(\af^*,0,\rho^*)$. 
\\

One must be careful at this point because the $k$-curve is supposed to reach the point $(\pi/2,\pi/2,\rho_1)$ at $k=1$. In fact, the curves $\CC^{\pm}$ may be considered as intersections of two surfaces $\CC^{\pm}_\af$ in $\R^3$ with parallel planes at level $\af$, $\af\in(0,\pi/2]$. By tracking back the point $\rho^*\in\CC^+\cap\CC^-$ at $\af=\af^*$, a space curve $\sfC\subset\CC^+_\af\cap\CC^-_\af$ is described. From [9, p. 147], the ``height'' $\af$ applied to $\sfC$ may be considered as a Morse function $f$. Therefore, whenever the tracking dies off, this represents a local extreme of $f$. Namely, the curve $\sfC$ is either descending to a smaller $\af$, or re-taking its ascent to a bigger $\af$. Moreover, $\sfC$ cannot be closed because $\sharp(\CC^+\cap\CC^-)$ is even for $\af<\af^*$ and $\CC^-$ shrinks to a point when $\af\to 0$.
\\

At the beginning of Section 5 we remarked that $\af$ could be restricted to the interval $(0,\pi/2]$. Letting $\af$ vary in $[\pi/2,\pi)$ we obtain the corresponding $\CC^{\pm}$ curves by reflection in $i\R$ of the $\CC^{\pm}$ curves for the $(0,\pi/2]$ case. In particular, $\CC^+\cap\CC^-=\emptyset$ when $\af$ approaches $\pi$. This will symmetrically extend $\sfC$ by $180^\circ$-rotation around the line $\af=\bt=\pi/2$, and make it connect $(\af^*,0,\rho^*)$ with $(\pi-\af^*,0,-\rho^*)$. Moreover, $(\pi/2,\pi/2,\rho_1)\in\sfC$ because its intersection number is odd. Namely, $(\af(k),\bt(k),\rho(k))$ is in fact a parametrisation of $\sfC$. This guarantees that we get a maximal {\it continuous} family $\SK_k$, $k\in(-\infty,\infty)$, passing through the standard Scherk-Karcher surface at $k=1$.    
\\

From this point on we shall strongly use the references [11] and [30]. Figures 32 and 33 of [30] are very helpful to understand what will happen to $\SK_k$. In [11] by the position of SGOH in $\R^3$, the Gau\ss \ map takes on the value $-1$ at the crossings between vertical and horizontal lines. These crossings correspond to points $A$ and $B$ of Figure 3(a). Therefore, we introduce the function
\[
   G=\frac{1+g}{1-g}.
\]      

Namely, $G$ is the same Gau\ss \ map, now clockwise rotated of $180^\circ$ about $Ox_2$. Because of that, the height differential must be taken as
\[
   dH=i\sqrt{\sin\af/2}\frac{(z-\ld)(z-\Ld)}{(z-\rho e^{i\bt})(z-\rho e^{-i\bt})}\cdot\frac{dz}{uz},
\] 
where the factor $\sqrt{\sin\af/2}$ is just a re-scaling and $\lambda^{\pm}$ are the $z$-values at which $g=1$, namely
\[
   [\kappa+2\cos\af\pm\sqrt{(\kappa+2\cos\af)^2-4}]/2. 
\] 

Since $\Lim{k\to\infty}{(\af,\bt,\rho)}=(\af^*,0,\rho^*)$, then $\Lim{k\to\infty}{\kappa}=\rho^*+1/\rho^*-2\cos\af^*$ and consequently $\Lim{k\to\infty}{(\ld,\Ld)}=(\rho^*,1/\rho^*)$. We recall that $T$ is the torus defined in Section 4. Denote its universal covering by $\mu:\C\to T$.
\\

In [11] one has some important functions and parameters that we shall use here by the same name, but in bold style to avoid confusion with our notation. So take $\bz$, $\bw$ and $\brho$ from [11] and notice that $\bz=iz$, $\brho=\pi/2-\af$, $2u^2=i\sin\af(\bz'/\bz)^2$ and $u^2\bw^2=2i\sin\af$. Without loss of generality, $u=e^{i\pi/4}\sqrt{\sin\af/2}\cdot\bz'/\bz$ and $u=-e^{i\pi/4}\sqrt{2\sin\af}/\bw$. Consequently,
$$
   G=\frac{\bw-\br e^{i\pi/4}}{\bw+\br e^{i\pi/4}}\leqno(29)
$$ 
and
$$
   dH=e^{i\pi/4}\frac{(\bz-i\ld)(\bz-i\Ld)}
   {(\bz-i\rho e^{i\bt})(\bz-i\rho e^{-i\bt})}\cdot\frac{d\bz}{\bz'},\leqno(30)
$$
where $\br=\sqrt{2\sin\af/\kappa}$.
\\

Now take any compact $\K\subset\C\setminus(\bz\circ\mu)^{-1}(i\rho^*)$. On $\K$, one sees that (29) and (30) will converge uniformly to the Weierstra\ss \ data of the SGOH, as presented in [11]. Since periods are closed for any real $k$, so they are for $k=\infty$. Moreover, the limit {\it must} be the SGOH from Hoffman-Karcher-Wei, since it is unique according to [8].
\\

For the embeddedness, standard arguments show that the fundamental piece of $\SK_1$ is a graph (see [15, p. 60], [27, p. 360] or [28, p. 566] for instance). Since the (maximal) family $\SK_k$ is continuous, a direct application of the {\it classical maximum principle} and {\it the maximum principle at infinity} shows that every member of this family {\it is} embedded, including the SGOH at the extremes. This completes assertion (c) in Theorem 1.1.
\\
\\
\\
\centerline{\bf 7. Appendix}
\\
\\
{\bf A.} Equality (22) holds if and only if (26) holds, for $\af=\bt=\pi/2$. But in this case (26) is equivalent to $\J_1=\J_2$, where
$$
   \J_1:=\frac{\kappa}{2}\int_0^{\pi/2}
   \frac{\sqrt{\cos t}\,dt}{\rho^2+1/\rho^2+2\cos(2t)}
   \eh\eh{\rm and}\eh\eh
   \J_2:=\int_0^{\pi/2}
   \frac{\sqrt{\cos t}^3dt}{\rho^2+1/\rho^2+2\cos(2t)},\leqno(31)
$$
with $\kappa=\rho-1/\rho$. Of course, if $\kappa\ge 2$ then $\J_1>\J_2$. Now observe that the $\rho$-derivative of $\J_1$ is positive when $\rho\in[1,1+\sqrt{2}]$, while $\J_2$ always decreases with $\rho$. Since $\rho\ge 1+\sqrt{2}$ implies $\kappa\ge 2$, then (31) can only hold for a {\it single} $\rho_1\in[1,1+\sqrt{2}]$. This is indeed the case, because in [11] the authors show that (31) has {\it at least} one solution.
\\
\\
{\bf B.} Since 
\[ 
   I_2=\int_0^\af\biggl(\frac{\cos\af-\cos\bt}{\cos t-\cos\af}\biggl)^\m\cdot\frac{dt}{\cos t-\cos\bt},
\]
the change of variables $u=\cos t$ gives 
\[ 
   I_2=(\cos\af-\cos\bt)^\m\int_{\cos\af}^1\frac{du/\sqrt{1-u^2}}{(u-\cos\af)^\m(u-\cos\bt)}.
\]

Since $1-u^2>1-u$ for $\af\in(0,\pi/2]$, the change $u=1-v^2$ leads to
\[
   \limsup_{\bt\to\af}I_2\le
   \limsup_{\bt\to\af}E\int_0^{(1-\cos\af)^\m}
   \frac{[(1-\cos\af)^\m-v]^{-\m}dv}{(1-\cos\bt)^\m-v},
\]
where $E:=2(\cos\af-\cos\bt)^\m(1-\cos\af)^{-\m}(1-\cos\bt)^{-\m}$. With $v=-w^2+(1-\cos\af)^\m$ we have
\[
   \limsup_{\bt\to\af}I_2\le
   \limsup_{\bt\to\af}2E\int_0^{(1-\cos\af)^{1/4}}
   \frac{dw}{(1-\cos\bt)^\m-(1-\cos\af)^\m+w^2}.
\]

This last integral gives $\limsup_{\bt\to\af}I_2\le2\pi\sqrt{2}/(1-\cos\af)^{3/4}$.  
\\
\\
{\bf C.} Now we are going to study the case $\rho\to\infty$. From (17) it follows that $\kappa\to\infty$, and so $q(t)$ vanishes uniformly on every compact subinterval of $[0,\af)$. However, $p(t)$ goes to infinity unless we re-scale $dh$ by $d\tilde{h}:=(1/\rho+\rho)dh$. Henceforth, we take this replacement for granted. The solution $w$ will continuously depend on the parameters at $t=\af$ providing the coefficients do {\it not} get singularities there. We can partially go round this by making the change $t(s)=\af+s^3$, $-\sqrt[3]{\af}<s<0$. Consequently, (9) becomes $\ddot{w}=P+Qw+R\dot{w}$, where $P(s)=9s^4p(t(s))$, $Q(s)=9s^4q(t(s))$ and $R(s)=3s^2r(t(s))+2/s$. Notice that $R$ still has singularities at $-\sqrt[3]{\af}$ and $0$, but this will not harm our analysis.\\

By fixing a compact subinterval $[-\sqrt[3]{\af},-\delta]$ of $[-\sqrt[3]{\af},0)$, $\rho\to\infty$ leads to $\ddot{w}=\tilde{p}+\tilde{r}\dot{w}$, where 
\[
   \tilde{p}(s)=\frac{-9s^4\sin 2t}{2(\cos t-\cos\af)};
\]
\[
   \tilde{r}(s)=3s^2\biggl(\cot t+\frac{0.5\sin t}{\cos t-\cos\af}\biggl)+\frac{2}{s}.
\]

Therefore,
\[
   \dot{w}(s)=\frac{-9s^2\sin t}{\sqrt{\cos t-\cos\af}}\int_{-\sqrt[3]{\af}}^s\frac{\sigma^2\cos t\,d\sigma}{\sqrt{\cos t-\cos\af}}.
\]

It is clear that $\dot{w}$ is always negative, and from $\ddot{w}=\tilde{p}+\tilde{r}\dot{w}$ the same holds for $\ddot{w}$. By using the continuous dependence on parameters, there exits $M>1$ such that both $\dot{w}$ and $\ddot{w}$ are still negative for any $\rho\ge M$. Since $Q$ is bounded at $s=0$, we can suppose that $M$ is big enough to keep $\ddot{w}<0$ for $s\in[-\delta,0]$, and consequently $\dot{w}(0)<0$. This means that $w'(\af)=-\infty$.     
\\

This should not surprise the careful reader, for although $w$ is finite on bounded curves, its derivatives can drastically change after re-parametrisation. Therefore, $\rho=1$ implies $w'(\af)>0$ while $\rho\to\infty$ implies $w'(\af)<0$. The {\it intermediate value theorem} assures the existence of a curve $\CC^+$ in the region $\A:=\{z\in\C:|z|>1>1-Im\{z\}\}$ with $w'|_{\CC^+}\equiv 0$. Since $w$ is real analytic, this curve separates $\A$ into a finite number of simply connected regions. The convergences are uniform for any $\bt\in[0,\pi]$, so that $\CC^+$ is bounded. 
\\

Similarly, for the stretch $t\in[0,\pi-\af]$ we take $t=\pi-\af+s^3$, $-\sqrt[3]{\pi-\af}<s<0$, so that $\rho\to\infty$ implies
\[
   \dot{w}(s)=\frac{-9s^2\sin t}{\sqrt{\cos t+\cos\af}}\int_{-\sqrt[3]{\pi-\af}}^s\frac{\sigma^2\cos t\,d\sigma}{\sqrt{\cos t+\cos\af}}.
\]

Now notice that $\dot{w}(s)$ varies from negative to positive for $s$ close to zero, when $\af$ drops from $\pi/2$ to zero. Anyway, $|\dot{w}(0)|=\infty$ except for a single $\af$-value in $(0,\pi/2)$.  
\\
\\
{\bf D.} Here we prove that $w'(\af)<0$ for $\bt=\pi$ and any $\rho\ge 1$. Due to the rotational symmetry around $X(\overline{AB})$, the ODE (9) for $z(t)=e^{it}$ and $0\le t\le\af$ can be equivalently studied for $-\af\le t\le 0$. In this case we have
\[
   w''-r(t)w'-q(t)w=p(t)\ge 0.
\]

By taking $s=-t/\af$, one rewrites the previous ordinary differential inequality as
\[
   \ddot{w}-R(s)\dot{w}-Q(s)w\ge 0,
\] 
where $Q(s)=\af^2q(t(s))$ and $R(s)=-\af r(t(s))$. Our intention is to apply Theorem 3.1. If $\bt=\pi$ and $\rho\ge 1$, then $\kappa=\rho+1/\rho+2\cos\af$ and 
\[
   \lim_{s\to 1}\biggl(-R+(1-s)Q\biggl)=-\infty.
\]

Moreover, the unitary normal is $(0,0,-1)$ at $z=e^{-i\af}$, and the third coordinate of $X(e^{-i\af})$ is negative. Therefore, $w$ has a positive maximum at $s=1$ and consequently $\dot{w}(1)>0$. By recalling that we took $-\af\le t\le 0$, then $w'(\af)<0$. 
\\
\\
\\
\centerline{\bf References}
\\
\begin{description}
\itemsep = 0.0 pc
\parsep  = 0.0 pc
\parskip = 0.0 pc
\item{[1]} \hspace{.04in} {\it F. Baginski}, Special functions on the sphere with applications to minimal surfaces, Special issue in honour of Rodica Simion, Adv. in Appl. Math. {\bf 28} (2002), 360--394.
\item{[2]} \hspace{.04in} {\it T. J. I'A. Bromwich}, A note on minimal surfaces, Proc. London Math. Soc. {\bf 30} (1899), 276--281.
\item{[3]} \hspace{.02in} {\it P. Collin}, Topologie et courbure des surfaces minimales propement plong\'ees de $\R^3$, Ann. of Math. {\bf 145} (1997), 1--31.
\item{[4]} \hspace{.02in} {\it C.J. Costa}, Example of a complete minimal surface in $\R^3$ of genus one and three embedded ends, Bol. Soc. Brasil. Mat. {\bf 15} (1984), 47--54.
\item{[5]} \hspace{.025in} {\it C.C. Chen} and {\it F. Gackstatter}, Elliptische und hyperelliptische Funktionen und vollst\"andige Minimalfl\"achen vom Enneperschen Type, Math. Ann. {\bf 259} (1982), 359--365.
\item{[6]} \hspace{.03in} {\it G. Darboux}, Le\c cons sur la th\'eorie g\'en\'erale de surfaces et les applications g\'eometri- ques au calcul infinit\'esimal, Gauthier-Villars, Paris, 2nd ed, 1914.
\item{[7]} \hspace{.01in} {\it J. Douglas}, Solution of the problem of Plateau, Trans. Amer. Math. Soc. {\bf 33} (1931), 263--321.
\item{[8]} \hspace{.025in} {\it L. Ferrer} and {\it F. Mart\'\i n}, Minimal surfaces with helicoidal ends, Math. Z. {\bf 250} (2005), 807--839.
\item{[9]} \hspace{.01in} {\it M.W. Hirsch}, Differential Topology, Springer, New York, 1976.
\item{[10]} {\it A. Huber}, On subharmonic functions and differential geometry in the large, Comment. Math. Helv. {\bf 32} (1957), 13--72.
\item{[11]} {\it D. Hoffman}, {\it H. Karcher} and {\it F. Wei}, The singly periodic genus-one helicoid, Comment. Math. Helv. {\bf 74} (1999), 248--279.
\item{[12]} {\it D. Hoffman}, {\it H. Karcher} and {\it F. Wei}, The genus one helicoid and the minimal surfaces that led to its discovery, Global analysis in modern mathematics, Publish or Perish, Houston, 119--170, 1993. 
\item{[13]} {\it D. Hoffman} and {\it W.H. Meeks}, Complete embedded minimal surfaces of finite total curvature, Bull. Amer. Math. Soc. {\bf 12} (1985), 134--136.
\item{[14]} {\it L. Hauswirth} and {\it M. Traizet}, The space of embedded doubly periodic minimal surfaces, Indiana Univ. Math. J. {\bf 51} (2002), no. 5, 1041--1079.
\item{[15]} {\it H. Karcher}, Construction of minimal surfaces, Surveys in Geometry, University of Tokyo 1--96, 1989 and Lecture Notes {\bf 12}, SFB256, Bonn, 1989.
\item{[16]} {\it H. Karcher}, The triply periodic minimal surfaces of Alan Schoen and their constant mean curvature companions, Manus- cripta Math. {\bf 64} (1989), 291--357.
\item{[17]} {\it H. Karcher}, Embedded minimal surfaces derived from Scherk's examples, Ma- nuscripta Math. {\bf 62} (1988), 83--114.
\item{[18]} {\it F.J. L\'opez} and {\it F. Mart\'\i n}, Complete minimal surfaces in $\R^3$, Publicacions Matematiques {\bf 43} (1999), 341--449.
\item{[19]} {\it F. Mart\'\i n} and {\it V. Ramos Batista}, The embedded singly periodic Scherk-Costa surfaces, Math. Ann. {\bf 336} (2006), no. 1, 155--189.
\item{[20]} {\it W.H. Meeks} and {\it H. Rosenberg}, The geometry of periodic minimal surfaces, Comment. Math. Helv. {\bf 68} (1993), 538--578.
\item{[21]} {\it J.C.C. Nitsche}, Lectures on minimal surfaces, Cambridge University Press, Cambridge, 1989.
\item{[22]} {\it R. Osserman}, A survey of minimal surfaces, Dover, New York, 2nd ed, 1986.
\item{[23]} {\it M.H. Protter} and {\it H.F. Weinberger}, Maximum principles in differential equations, Prentice-Hall, 1967.
\item{[24]} {\it T. Rad\'o}, On the problem of Plateau, Ergeben. d. Math. u. ihrer Grenzgebiete, Springer-Verlag, Berlin, 1933.
\item{[25]} {\it V. Ramos Batista}, Singly periodic Costa surfaces, J. London Math. Soc. {\bf 72} (2005), no. 2, 478--496. 
\item{[26]} {\it V. Ramos Batista}, Noncongruent minimal surfaces with the same symmetries and conformal structure, Tohoku Math. J. {\bf 56} (2004), 237--254.
\item{[27]} {\it V. Ramos Batista}, A family of triply periodic Costa surfaces, Pacific J. Math. {\bf 212} (2003), 347--370.
\item{[28]} {\it V. Ramos Batista}, The doubly periodic Costa surfaces. Math. Z. {\bf 240} (2002), 549--577.
\item{[29]} {\it H.W. Richmond}, On minimal surfaces, Jour. London Math. Soc. {\bf 19} (1944), 229--241.
\item{[30]} {\it M. Weber}, The genus one helicoid is embedded, Habilitation Thesis, Bonn, 2000. 
\end{description}
\centerline{\underline{\hspace{3cm}}}
\bigskip
\bigskip
\centerline{Department of Mathematics, the George Washington University} 
\centerline{ Hall of Government 224, 2115 G Str. NW, Washington DC 20052, USA}
\centerline{e-mail: baginski@gwu.edu}
\bigskip
\centerline{CMCC - UFABC, r. Catequese 242, 09090-400 Santo Andr\'e-SP, Brazil}
\centerline{e-mail: valerio.batista@ufabc.edu.br}
\end{document}